\documentclass[11pt,leqno]{article}
\usepackage{amsthm,amsmath,amsfonts,latexsym,amssymb}
\usepackage{mathrsfs,dutchcal,arydshln,gensymb,appendix,abstract}
\usepackage{graphicx}
\usepackage{tikz}
\date{}

\begin{document}
	\title{A family of area-minimizing tensor varieties}
	\author{\centerline {Hongbin Cui}}
 
 \maketitle
\begin{center}
        \itshape \small Dedicated to Anguo Cui
    \end{center}
 \begin{abstract}
By using Lawlor's curvature criterion, we prove that a class of real tensor varieties—tensors of \textit{multilinear rank} at most $(1, \ldots, 1)$—are regular area-minimizing cones except for one case, which gives a partial generalization of their minimality as in \cite{HKV23}. Such area-minimizing tensor varieties have not been found beyond matrix varieties (cf. \cite{KL99} and \cite{CJX24}). Moreover, we provide a new derivation for the important Lawlor ODE.
 \end{abstract}

		
		
		\bigskip
		
		\tableofcontents

\section{Introduction}
\subsection{Backgrounds on area-minimizing cones}

\medskip\noindent

Area-minimizing surfaces in $\mathbb{R}^{n}$ are surfaces (or integral currents, cf. \cite{FF60}) that globally minimize the area functional with given boundaries. As solutions to Plateau's problem, these surfaces often exhibit singularities in the interior. The celebrated regularity theorems provide quantitative descriptions for the Hausdorff dimension of the interior singular set (see a nice survey in \cite{M16}). However, understanding the concrete structures and 'topological type' of the singularities is of equal importance. It is well-known that, at every interior singular point of an area-minimizing surface, there always exists at least one tangent cone as the \textit{blowing-up} limit, which can be seen as the first-order behavior of singularities—just as tangent planes characterize the local behavior of smooth manifolds. This cone is itself area-minimizing (cf. Sections 4.3.16 and 5.4.3 in \cite{Fe69}). Therefore, studying and finding area-minimizing cones help understand the singularities of area-minimizing surfaces. Moreover, many new examples of area-minimizing hypersurfaces can be produced by perturbing regular area-minimizing hypercones (cf. \cite{HS85}). We call a cone $\mathcal{C} \subset \mathbb{R}^{n}$ a \textit{regular cone} if it intersects $\mathbb{S}^{n-1}$ in a closed smooth submanifold, so the origin is the only singularity. This concept was first defined in \cite{HS85} for cones over minimal compact smooth hypersurfaces in the sphere $\mathbb{S}^{n-1}$. And any other cones with larger singular sets will be simply called \textit{non-regular cones}. However, at least in the cases of area-minimizing hypercones, one can repeatedly apply the "blow-up" procedures to handle singularities away from the origin, and finally derive regular area-minimizing cones (cf. \cite[Proposition 9.6]{G84}) in $\mathbb{R}^{n-k}$ for some integer $k$. 

\medskip\noindent

Area-minimizing cones also play important roles in the fields of differential geometry and partial differential equations. The celebrated Bernstein problem asks whether an entire solution to the minimal surface equation (MSE) defined on $\mathbb{R}^n$ must be flat. It seemed impossible to obtain satisfactory answers to this question using only the methods of partial differential equations or complex analysis in the case $n=2$ before, since there exists a critical dimension: $n=7$. The research on area-minimizing cones completely solved this question (cf. a discussion in \cite[p.102-103]{CM11}): Fleming (\cite{Fl62}) found that the \textit{blow-down} of a non-flat entire solution to MSE on $\mathbb{R}^n$ implies a non-flat area-minimizing cone in $\mathbb{R}^{n+1}$, which was then shown by De Giorgi \cite{DG65} to be of cylindrical type $\mathcal{C}\times \mathbb{R}$ (also see \cite[Proposition 9.6]{G84}), thus reducing the Bernstein problem to the study of area-minimizing cones in $\mathbb{R}^n$.
In 1968, J. Simons (\cite{S68}) proved that no stable minimal hypercones exist in $\mathbb{R}^n$ for $n\leq 7$, proving affirmatively in the cases: $n\leq 7$. In 1969, Bombieri, De Giorgi and Giusti (\cite{BDGG69}) found the first non-flat regular area-minimizing cone $\mathcal{C}(S^3 \times S^3)$ (\textit{Simons cone}) in $\mathbb{R}^8$, thus finally settled the Bernstein problem. 

\medskip\noindent

In high-codimension cases, Y. Yuan (\cite{Y02}) used similar tangent cone analysis to prove an associated Bernstein theorem for special Lagrangian equations. Moreover, high-codimension area-minimizing cones play important roles in the study of Dirichlet problems for minimal surface systems (cf. \cite{LO77}, \cite{XYZ19}). Additionally, minimal cones are natural critical points for general radially weighted area functionals, and area-minimizing cones have also been shown to be minimizers of several such weighted functionals, such as the $\textbf{E}$-functional in mean curvature flow and the Euler-Dierkes-Huisken functional (a weighted area functional introduced in \cite{DH24}), see \cite{D20}, \cite{DH24},\cite{CX25}.

\medskip\noindent

Due to these important properties and applications, studying and finding area-minimizing cones (regular or non-regular) has attracted considerable attention as in \cite{La72}, \cite{Si74}, \cite{HL82}, \cite{FK85}, \cite{Lin87}, \cite{M83}, \cite{Ma87}, \cite{La89}, \cite{La91}, \cite{W94}, \cite{La98}, \cite{KL99}, \cite{Da04}, \cite{DP09}, \cite{Z16}, \cite{Z18} and references therein. 
		
\medskip\noindent 

A general method for proving the area-minimization of a given regular minimal cone, called the \textit{curvature criterion}, was provided by Lawlor (cf. \cite{La91}). See Sect. 2. His method is special in that the calculations can only be carried out within the framework of differential geometry. These calculations involve the second fundamental forms of the smooth links in unit spheres and a specific type of 'embedded radius' of the cones. In most cases, the method ultimately reduces to calculating two important quantities: the \textit{vanishing angle} and \textit{normal radius}. As an application, Lawlor completely classified the area minimizations for cones over arbitrary products of spheres. Since then, regular area-minimizing cones of various topological types can be found from the research of minimal submanifolds in spheres, as in \cite{Ke94}, \cite{Ka02}, \cite{XYZ18}, \cite{TZ20}, \cite{OS21}, \cite{JC22}, \cite{JCX22}, \cite{CJL25}, etc. 
		
\medskip\noindent 

For non-regular area-minimizing cones, projective varieties provide many such examples when considered as cones in ambient Euclidean spaces. Via calibrated geometry theory (\cite{HL82}), complex algebraic varieties in ambient spaces $\mathbb{C}^{n} \cong \mathbb{R}^{2n}$ are area-minimizing in any of their compact portions. However, for real projective varieties, the situation becomes complicated: firstly, real projective varieties are not necessarily minimal surfaces. Secondly, they are often non-regular cones, and the Lawlor curvature criterion does not apply. In 1998, Lawlor \cite{La98} developed another powerful but more general method, called \emph{directed slicing}, to study highly singular minimizing surfaces. 
\medskip\noindent

\subsection{Real matrix varieties and tensor varieties}

\medskip\noindent 

In \cite{KL99}, by using the method of directed slicing, a family of matrix varieties—the real determinantal varieties---were proved to be area-minimizing except for two cases. Recently, similar results for Pfaffian varieties were also established in \cite{CJX24}.

\medskip\noindent 

It is natural to consider some suitable generalizations to real tensor varieties from these real matrix varieties. Let $\mathbb{R}^{n_1} \otimes \cdots \otimes \mathbb{R}^{n_d}$ be the space of real $\left(n_1, \ldots, n_d\right)$-tensors. By endowing each $\mathbb{R}^{n_j}$ with the orthonormal bases $\mathcal{E}^j=\left\{\mathbf{e}_1^j, \ldots, \mathbf{e}_{n_j}^j\right\}$, every tensor $T$ can be written as:
$T=\sum_{i_1=1}^{n_1} \cdots \sum_{i_d=1}^{n_d} t_{i_1 \cdots i_d} \mathbf{e}_{i_1}^1 \otimes \cdots \otimes \mathbf{e}_{i_d}^d$. $\mathbb{R}^{n_1} \otimes \cdots \otimes \mathbb{R}^{n_d}$ is equipped with the Frobenius inner product: $
\langle T, S\rangle=\sum_{j=1}^d \sum_{i_j=1}^{n_j} t_{i_1 \ldots i_d} s_{i_1 \ldots i_d}
$
where $T=\left(t_{i_1 \cdots i_d}\right)$ and $S=\left(s_{i_1 \cdots i_d}\right)$. The rank of matrices and the determinantal varieties have the following generalizations to tensors:

\medskip\noindent 

\textbf{Definition 1.1. (\cite{HKV23}, \cite[p.33, 35, p.69 and p.102]{Lan12})} \textit{Define the multilinear rank (or mode $j$ rank) of $T \in \mathbb{R}^{n_1} \otimes \cdots \otimes \mathbb{R}^{n_d}$ to be the $d$-tuple of numbers:
$$
\mathbf{R}_{\operatorname{multlin}}(T):=\min \left\{\left(r_1, \ldots, r_d\right) \in \mathbb{N}_0^d: T \in W_1 \otimes \cdots \otimes W_d, W_j \subseteq \mathbb{R}^{n_j}, \operatorname{dim} W_j=r_j\right\}
$$
where $\mathbb{N}_0$ denotes the set of positive integers and each $W_j$ is a linear subspace of $\mathbb{R}^{n_j}$. The subspace varieties $\hat{S} u b_{r_1, \ldots, r_d}\left(\mathbb{R}^{n_1} \otimes \cdots \otimes \mathbb{R}^{n_d}\right)$ is  defined as
\begin{equation}\notag
    \hat{S} u b_{r_1, \ldots, r_d}\left(\mathbb{R}^{n_1} \otimes \cdots \otimes \mathbb{R}^{n_d}\right): 
    =\left\{T \in \mathbb{R}^{n_1} \otimes \cdots \otimes \mathbb{R}^{n_d} \mid \mathbf{R}_{\operatorname{multlin}}(T) \leq\left(r_1, \ldots, r_d\right)\right\}.
\end{equation}
Also, define the rank of a tensor $T$:
$$
\textbf{R}(T):=\min \left\{r \in \mathbb{N}_0: T=\sum_{i=1}^r v_{1} \otimes \cdots \otimes v_{d}, v_{j} \in \mathbb{R}^{n_j}\right\},
$$
then the space of rank one tensors is just the subspace variety with $\left(r_1, \ldots, r_d\right)=(1,\ldots,1)$ up to a zero element, which is called the $d$-factor Segre variety, abbreviated as Segre variety here.}

\medskip\noindent 

The number $r_j$ is equal to the rank of the $n_j \times \prod_{k \neq j} n_k$ matrix obtained by flattening or unfolding the tensor along mode $j$ (cf. \cite[Theorem 6.13]{H12}). Also, it satisfies (cf. \cite[p. 36]{Lan12})
$$
 r_j \leq \prod_{k \neq j} r_k.
$$

\medskip\noindent 

Similar to determinantal varieties and Pfaffian varieties, the subspace varieties are stratified by the multilinear rank. Denote its regular part by $\mathcal{T}_{d, \mathbf{r}}$, i.e. the space of tensors of multilinear rank $\mathbf{r}=\left(r_1, \ldots, r_d\right)$. In \cite{HKV23}, Heaton, A., Kozhasov, K. and Venturello, L. proved that subspace varieties are minimal varieties:

\medskip\noindent

\textbf{Proposition 1.2. (\cite{HKV23})} \textit{$\mathcal{T}_{d, \mathbf{r}}$ are smooth minimal submanifolds in $\mathbb{R}^{n_1} \otimes \cdots \otimes \mathbb{R}^{n_d}$.}

\medskip\noindent 

\textbf{Remark 1.3.} \textit{In the notations of Definition 1.1, determinantal varieties: $\hat{S} u b_{r, r}\left(\mathbb{R}^{n_1} \otimes \mathbb{R}^{n_2}\right)$ and Segre varieties: $\hat{S} u b_{1, \ldots, 1}\left(\mathbb{R}^{n_1} \otimes \cdots \otimes \mathbb{R}^{n_d}\right)$ are the two simplest cases of subspace varieties. The area-minimization of determinantal varieties was studied in \cite{P93} and \cite{KL99}. It is reasonable to investigate the area-minimization of general subspace varieties using the direct slicing method. However, the problem is challenging for several reasons: the \textit{standard forms} and the low-rank approximation theory of matrices used in \cite{KL99} and \cite{CJX24} are not well-established for tensors.}

\medskip\noindent

\subsection{Statement of the main results}

\medskip\noindent

In this paper, we study Segre varieties. Unlike determinantal varieties $\hat{S} u b_{r, r}\left(\mathbb{R}^{n_1} \otimes \mathbb{R}^{n_2}\right)$, Segre varieties are regular cones whose regular parts are $\mathcal{T}_{d, (1,\cdots,1)}$. Thus, we can apply Lawlor's curvature criterion. Using the local parametrization given in Sect. 3, Segre varieties can be written as the cones: 
$$
\mathcal{C}=C(\mathbb{S}^{n_1-1} \otimes \cdots \otimes \mathbb{S}^{n_d-1})=\{0\} \cup \mathcal{T}_{d, (1,\cdots,1)}.
$$

\medskip\noindent

And we prove that

\medskip\noindent 
 
 \textbf{Main Theorem}. {\it Except for $C(\mathbb{S}^{1} \otimes \mathbb{S}^{1} \otimes \mathbb{S}^{1})$, all other Segre varieties $\mathcal{C}=C(\mathbb{S}^{n_1-1} \otimes \cdots \otimes \mathbb{S}^{n_d-1})$ are regular area-minimizing cones.}
 
\medskip\noindent 

For the case  $C(\mathbb{S}^{1} \otimes \mathbb{S}^{1} \otimes \mathbb{S}^{1})$, the vanishing angle does not exist, and Lawlor's curvature criterion does not apply.

\medskip\noindent 

The organization of the paper is as follows. In Sect. 2, we review and discuss Lawlor's curvature criterion and give a new derivation for the important Lawlor ODE through a reverse process of Lawlor's. The inequality for $J_{k}\Pi$, which seems to have been overlooked in Lawlor's original treatment, is accurately derived. In Sect. 3, we introduce the local parameterization of the Segre variety and calculate the maximal square norm of the shape operator and the normal radius. The proof of \textit{Main Theorem} then follows.

\medskip\noindent
\medskip\noindent

\section{Lawlor's curvature criterion}

\medskip\noindent 

Let $\mathcal{C}=C(\Sigma) \subset \mathbb{R}^{n}$ be a $k$-dimensional regular minimal cone with isolated singularity at the origin, where $\Sigma $ is a smooth, minimal submanifold in $S^{n-1}(1)$ and $k \geq 3$. Firstly, we recall Lawlor's definitions for \textit{normal wedge} and \textit{normal radius}:

\medskip\noindent 

\textbf{Definition 2.1. (\cite[Definition 1.1.1,1.1.2,1.1.3]{La91})} \textit{\emph{(1)} For $p\in \Sigma$ and $\alpha \geq 0$, let $U_{p}(\alpha)$ be the union of all open normal geodesics (in the sphere) of length $\alpha$ from $p$. The \textit{normal wedge} $W_{p}(\alpha)$ is defined as the conical set over $U_{p}(\alpha)$.}

\textit{\emph{(2)} Let $N_{p}:= {\rm max} \{\alpha | W_{p}(\alpha)\cap \mathcal{C}= \overrightarrow{op}\}$ , where $\overrightarrow{op} \in  \mathcal{C}$ is the ray through $p\in \Sigma$, we call $N_{p}$ the \textit{normal radius} at $p$.}

\medskip\noindent

Choose a unit normal vector $v$ at $p$, denote the $\mathbf{2}$-dimensional subspace spanned by $\{\overrightarrow{op}, v\}$ of $\mathbb{R}^{n}$ by $\mathcal{L}_{p,v}$, then:

\medskip\noindent

\textbf{Lemma 2.2.} \textit{Assume $q \in \mathcal{L}_{p,v} \cap \Sigma$ for some unit normal $v$, $q\neq p$, then $N_{p}$ is the infimum of the angle between $\overrightarrow{op}$ and $\overrightarrow{oq}$ among all chosen normals $v$. Equivalently, $N_{p}$ is the shortest normal geodesic in the unit sphere starting from $p$ that intersects another point of $\Sigma$.} 

\medskip\noindent

The normal radius controls the size of the angle neighborhood for which the following \textit{Lawlor retraction map} can be defined in: 

\medskip\noindent

%The normal exponential map of $\Sigma$ in the unit sphere is: ${\rm exp}(p,v)=\cos  |v| \ \overrightarrow{op}+ \sin  |v| \ \frac{v}{|v|}$, Lawlor consider an expansion of ${\rm exp}_{p}(\theta v)$ in the radial direction: $ S(p,v)=r(|v|) {\rm exp}(p,v)$, then extend $S(p,v)$ homothetically when $p$ is moving along the radial direction of $\mathcal{C}$.  

%\medskip\noindent

On $\mathcal{L}_{p,v}$, select a curve in polar coordinate $(r(\theta),\theta)$: $S_{p}=r(\theta) ( \cos  \theta \ \overrightarrow{op}+ \sin  \theta \ v)$, then extend $S_{p}$ homothetically to a map $\Phi$ from the cone $\mathcal{C}$ to the ambient Euclidean space. Lawlor wants to find some function $r(\theta)$ such that $S_{p}$ can extend to infinity when $\theta$ approaches to some fixed angle $\theta_{0}$ and also the $k$-dimensional Jacobian: $J_{k}\Phi \geq 1$, it implies that the resulting map $\Phi$ is area non-decreasing and $\Phi$ takes values at a normal wedge of $\mathcal{C}$ at $p$. And, if pointwise the normal wedges $ W_{p}(\theta_{0}(p))$ are disjoint for any $p,q \in \Sigma, p\neq q$:
\begin{equation}\label{Normal radius}
 W_{p}(\theta_{0}(p))\cap W_{q}(\theta_{0}(q))=\varnothing, 
\end{equation}
the inverse map $\Pi:=\Phi^{-1}$ is well-defined in an angle neighborhood of the cone $W:= \bigcup_{p\in \Sigma} W_{p}(\theta_{0}(p))$, then define that $\Pi$ maps every point outside $W$ to the origin, it follows that $\Pi$ is area non-increasing. 

\medskip\noindent

\textbf{Definition 2.3. (\cite[Sect. 1.1]{La91})} \textit{We call the above area non-increasing retraction map $\Pi$: the \textit{Lawlor retraction map}.}

\medskip\noindent

The existence of Lawlor retraction map (\cite[p.10-12]{La91}) is reduced to the following ordinary differential equation which we just call it \textit{Lawlor ODE}.

\medskip\noindent

\textbf{Theorem 2.4.} \textit{If
\begin{equation}\label{va1}
    \frac{dr}{d\theta}\leq r \sqrt{r^{2k}({\rm cos}\theta)^{2k-2}\left(\mathrm{inf}_{v}{\rm det}(I-{\rm tan} \theta A_v)\right)^{2}-1}
\end{equation}
for any normal direction $v$, then $J_{k}\Pi \leq 1$, where $A_v$ is the shape operator of $\Sigma \subset S^{n-1}$ in normal direction $v$.}

\medskip\noindent

Here we give a new deduction for Lawlor ODE through a reverse process of Lawlor's (\cite[Chapter 1]{La91}), which has its own significance: the inequality for $J_{k}\Pi$ is accurately given here and the concrete conditions for $J_{k}\Pi=1$ is also discussed. 

\medskip\noindent

\textbf{Proof of Theorem 2.4:} By the homothetical property of $\Phi$, we can only compute the Jacobian at a point $p\in \Sigma$. We make the most general assumption for $\Phi$: except for  $r=r (\theta)$ only, in the radial direction of the cone, we assume the angle $\theta$ is affected by radial movements (which is indeed the case!) However, it can be assumed that the normal vector $v$ remains unchanged in the radial direction. In the spherical tangential directions, we assume that $\theta $ and the extension of $v$ can be both affected by the movements of $p$.

\medskip\noindent

Denote the position vector by $x$, and let $\left\{p=v^0, v^1, \cdots, v^{k-1}\right\}$ denote  an orthonormal tangential basis of $C$ at $p$, we always identify the following: $p=v^0=\overrightarrow{o p}$, choose a ray starting from $p$ in the direction $p$: $\alpha(t)=p+t p$, then we have $\Phi(\alpha(t))=(1+t)(r(\theta) \cos \theta p+r(\theta) \sin \theta v)$, derivative in $t$, we get:
$$
\Phi_*\left(v^0\right)=\left(r \cos (\theta)+v^0(r \cos \theta)\right) p+\left(r \sin (\theta)+v^0(r \sin \theta)\right) v:=A p+B v,
$$
denote $a=v^0(\theta)=\frac{d \theta}{d t}(0)$, then
$$
A=\left(r \cos \theta+\left(\frac{d r}{d \theta} \cos \theta-r \sin \theta\right) a\right) {\rm and} \ B=\left(r \sin \theta+\left(\frac{d r}{d \theta} \sin \theta+r \cos \theta\right) a\right),
$$
it implies that
\begin{equation}\label{rjk}
\begin{aligned}
\left|\Phi_*\left(v^0\right)\right| & =\sqrt{A^2+B^2} \\
& =\left(\left(\frac{d r}{d \theta}\right)^2+r^2\right) a^2+2 r \frac{d r}{d \theta} a+r^2 \\
& \geq \frac{r^2}{\sqrt{\left(\frac{d r}{d \theta}\right)^2+r^2}},
\end{aligned}
\end{equation}
and  the minimum is attained iff $a=v^0(\theta)=\frac{-r \frac{d r}{d \theta}}{r^2+\left(\frac{d r}{d \theta}\right)^2}$, that recovers Lawlor's formula in \cite[p.9]{La91}.

\medskip\noindent

Next, we calculate the tangent map of $ \Phi$ in the spherical tangential directions. Assume $p$ moves on the minimal submanifold $\Sigma \subset S^{n-1}(1)$ along a curve $\sigma(t) \in \Sigma$, which satisfies $\sigma(0)=p, \sigma^{\prime}(0)=v^i$. We assume that there is an extension of $v$ along $\sigma(t)$: $\tilde{V}(\sigma(t)), \tilde{V}(p)=v$, then
$$
\begin{aligned}
\Phi_*\left(v^i\right) & =\left.\frac{d}{d t}\right|_{t=0}(r \cos \theta \sigma(t)+r \sin \theta \tilde{V}(\sigma(t))) \\
& =r \cos \theta v^i+r \sin \theta D_{v^i} \tilde{V}+v^i(r \cos \theta) p+v^i(r \sin \theta) v,
\end{aligned}
$$
where $D$ is Euclidean connection, also denote the Riemannian connection in $S^{n-1}(1)$ by $\nabla$, we can further decompose:
$$
\begin{aligned}
D_{v^i} \tilde{V} & =\nabla_{v^i} \tilde{V}+\left\langle D_{v^i} \tilde{V}, x\right\rangle x \\
& =\nabla_{v^i} \tilde{V}-\left\langle v, D_{v^i} x\right\rangle x \\
& =\nabla_{v^i} \tilde{V} \\
& =-A_v v^i+\nabla_{v^i}^{\perp} \tilde{V},
\end{aligned}
$$
where $A$ is the shape operator---independent with the extension of the normal vector $v$, $\nabla^{\perp}$ denotes the normal connection of $\Sigma$ in the unit sphere, it follows that
$$
\Phi_*\left(v^i\right)=r \cos \theta\left(I-\tan \theta A_v\right) v^i+r \sin \theta \nabla_{v^i}^{\perp} \tilde{V}+v^i(r \cos \theta) p+v^i(r \sin \theta) v.
$$

\medskip\noindent

Now, we choose an orthonormal normal vectors at $p$: $\left\{v=v^k, \cdots, v^{n-1}\right\}$, extending to a local normal frame: $\tilde{V}^\alpha(\alpha=k, \cdots, n-1)$. Let $F$ denote the $k$-vector $\Phi_*\left(v^0\right) \wedge \Phi_*\left(v^1\right) \wedge \cdots \wedge \Phi_*\left(v^{k-1}\right)$, then
$$
F =\left(A v^0+B v^k\right) \wedge \cdots \wedge\left(r \cos \theta\left(I-\tan \theta A_v\right) v^i+v^i(r \cos \theta) v^0+v^i(r \sin \theta) v^k +r \sin \theta \nabla_{v^i}^{\perp} \tilde{V}^k\right) \wedge \cdots,
$$
for convenience, we write $r \operatorname{cos} \theta\left(I-\tan \theta A_v\right) v^i=L_{i j} v^j(1 \leq i, j \leq k-1), m_i=v^i(r \cos \theta)$, $n_i=v^i(r \sin \theta)$. Since $\bigwedge^k\left(T_p \mathbb{R}^n\right)$ has an ordered orthonormal basis: $\left\{v^{s_0} \wedge \cdots v^{s_{k-1}} \mid 0 \leq s_0<\cdots<\right.$ $\left.s_{k-1} \leq n-1\right\}$. It follows that,

\medskip\noindent

$$
\begin{aligned}
F & =A \operatorname{det}(L) v^0 \wedge v^1 \wedge \cdots v^{k-1}+B \operatorname{det}(L) v^k \wedge v^1 \wedge \cdots v^{k-1} + \\
& \sum_{i=1}^{k-1}\left|\begin{array}{cccccc}
L_{11} & \cdots & \widehat{L_{1 i}} & \cdots & L_{1 k-1} & A n_1-B m_1 \\
\vdots & \cdots & \vdots & \cdots & \vdots & \vdots \\
L_{k-1,1} & \cdots & \widehat{L_{k-1, i}} & \cdots & L_{k-1, k-1} & A n_{k-1}-B m_{k-1}
\end{array}\right| v^0 \wedge \cdots \wedge \widehat{v^i} \wedge \cdots \wedge v^{k-1} \wedge v^k \\
& + \text {other terms},
\end{aligned}
$$
by Cramer's rule, the third terms are all zeros iff
$$
A n_i-B m_i=0,
$$
for any $1 \leq i \leq k-1$ which is further equivalent to $v^i(\theta)=0$. The \textit{other terms} depend on  $\nabla_{v^i}^{\perp} \tilde{V}^k$ and thus depend on how the normal vectors are extended, a sufficient condition for the vanishing of these terms is that $\nabla^{\perp} \tilde{V}^k=0$, i.e. $\nabla^{\perp} \tilde{V}=0$.

\medskip\noindent

Thus, together with \eqref{rjk}, we can conclude that 
\begin{equation}\label{Jk}
J_{k}\Phi = {\rm sup} \left|\Phi_*\left(v^0\right) \wedge \Phi_*\left(v^1\right) \wedge \cdots \wedge \Phi_*\left(v^{k-1}\right)\right| \geq \frac{r^2(r \cos \theta)^{k-1}}{\sqrt{\left(\frac{d r}{d \theta}\right)^2+r^2}} \operatorname{det}\left(I-\tan \theta A_v\right),
\end{equation}
it implies that 
$$
J_{k}\Pi \leq \frac{\sqrt{(d r / d \theta)^2+r^2}}{r^{k+1} (\cos \theta)^{k-1} \operatorname{det}\left(I-\tan \theta A_v \right)},
$$
the condition that the right hand is less than or equal to $1$ is equivalent to 
$$
\frac{dr}{d\theta}\leq r \sqrt{r^{2k}({\rm cos}\theta)^{2k-2}\left(\mathrm{inf}_{v}{\rm det}(I-{\rm tan} \theta A_v)\right)^{2}-1}.
$$
That is just Lawlor ODE in the equal sign case. Moreover, under the conditions: (1) the extended normal vector field $
\tilde{V}$ is normal parallel; 
(2) $r=r(\theta)$ solves the Lawlor ODE
and (3)
$$
v^0(\theta)=\frac{-r \frac{d r}{d \theta}}{r^2+\left(\frac{d r}{d \theta}\right)^2}, v^i(\theta)=0 \ (1 \leq i \leq k-1),
$$
$\Phi$ and $\Pi$ are area-preserving maps in a $\theta_{0}$-angle neighborhood of $\mathcal{C}$. $\Box$

\medskip\noindent

The inequality in \eqref{va1} can be replaced by the equality to make the possibility for the non-intersection of $\theta_0$-normal wedges bigger. Then, the solution is analyzed by Lawlor:  either $\frac{dr}{d\theta}$ vanishes at some positive $\theta(p)$ or $r\rightarrow \infty$ as $\theta \rightarrow \theta_{0}(p)$, in the latter case, Lawlor can construct his retraction map $\Pi$ if in addition \eqref{Normal radius} is satisfied for the normal wedges with the angle $\theta_{0}$, finally it follows that the cone is area-minimizing (see \cite[Theorem 1.2.1]{La91}).

\medskip\noindent
 
\textbf{Definition 2.5. (\cite[Definition 1.1.7]{La91})} \textit{We call the smallest $\theta_{0}(p)$ among all the normal directions $v$ the vanishing angle at $p$.}

\medskip\noindent

The global condition \eqref{Normal radius}  depends on the submanifold geometry of $\Sigma$ and is not easy to check directly. In fact, if $\Sigma$ is a minimal isoparametric hypersurface of the unit sphere, two normal wedges whose radii are smaller than the focal radius will never intersect, since the parallel hypersurfaces of $\Sigma$ are also embedded isoparametric hypersurfaces if the distance is less than the focal distance. Also, the positivity of ${\rm det}(I-{\rm tan} \theta A_v)$ implies that the resulting vanishing angles are always less than the focal angle. Thus, for isoparametric cases, we can only consider the existence of vanishing angles. For general cases, Lawlor gave the following simplified criterion related to the normal radius:

\medskip\noindent

\textbf{Theorem 2.6. (\cite[Theorem 1.3.5]{La91})} \textit{If the vanishing angle $\theta_{0}(p)$ exists for any points $p\in \Sigma$ and pointwisely satisfy
\begin{equation}\label{two times vanishing angle}
2\theta_{0}(p)\leq N(p),
\end{equation}
where $N(p)$ is the normal radius at $p$, then the cone $\mathcal{C}=C(\Sigma)$ is area-minimizing (in the sense of mod $2$ when $\Sigma$ is nonorientable).}

\medskip\noindent

\textbf{Remark 2.7.} \textit{In fact, we often use the estimated vanishing angle in the curvature criterion, which had been given in Lawlor's Table (cf. \cite[sect. 1.4]{La91}), the estimated vanishing angle only depend on the maximum values of the norm of shape operator $||A_{v}||^2$ and the dimension of the cone, see \cite[Sect. 3.1]{TZ20} for more information.}

\medskip\noindent
\medskip\noindent

\section{Area-minimizing Segre varieties}

\medskip\noindent

\textbf{3.1. The tangent spaces and normal spaces:} the tangent space of Segre variety is derived in \cite{HKV23} from the pointview of group action. However, we can easily find them by local parameterization. The local parameterization for the Segre varieties is given as follows (also see \cite[p.69]{Lan12}):
\begin{equation}\label{pa}
\begin{aligned}
f: [0,\infty) \times \mathbb{S}^{n_1-1} \times \cdots \times \mathbb{S}^{n_d-1} & \rightarrow C(\mathbb{S}^{n_1-1} \otimes \cdots \otimes \mathbb{S}^{n_d-1}) \\
\left(\lambda, e_1, \cdots, e_d\right) & \mapsto \lambda e_1 \otimes \cdots \otimes e_d,
\end{aligned}
\end{equation}
obviously $(\lambda,\cdots,-e_{i_{1}}, \cdots,-e_{i_{2k}},\cdots)$ maps to the same point.

\medskip\noindent

Given every $\mathbb{R}^{n_j}$ the orthonormal bases $\mathcal{E}^j=\left\{\mathbf{e}_1^j, \ldots, \mathbf{e}_{n_j}^j\right\}$, we choose the point $p=e_{1}^1 \otimes \cdots \otimes e_1^{d}$, by chain rules, choose a curve $e_{i}(t)$ ranges in $\mathbb{S}^{n_i-1}$ then derivative, we can see the tangent basis of $\mathcal{C}=C(\mathbb{S}^{n_1-1} \otimes \cdots \otimes \mathbb{S}^{n_d-1})$ at $p$ can be given as:
$$
T_{p}\mathcal{C}={\rm span}\{e_{1}^1 \otimes \cdots \otimes e_1^{d}, \cdots, e_{1}^1 \otimes \cdots \otimes e_{1}^{i-1} \otimes e_{\alpha_{i}}^{i} \otimes e_{1}^{i+1}  \otimes \cdots \otimes e_1^{d}, \cdots \},
$$
where $1 \leq i\leq d, 2\leq \alpha_{i} \leq n_{i} $.

\medskip\noindent

Denote $E_{\alpha_{i}}^{i}=e_{1}^1 \otimes \cdots \otimes e_{1}^{i-1} \otimes e_{\alpha_{i}}^{i} \otimes e_{1}^{i+1}  \otimes \cdots \otimes e_1^{d}$, $E_{\alpha_{i}\beta_{j}}^{ij}=e_1^1 \otimes \cdots \otimes e_{\alpha_{i}}^{i} \otimes \cdots \otimes e_{\beta_{j}}^{j} \otimes \cdots \otimes e_1^d: 2\leq \alpha_{i} \leq n_{i}, 2\leq \beta_{j} \leq n_{j}$ and general $E_{\alpha_{i_1}\cdots \alpha_{i_l}}^{i_{1}\cdots i_{l}}=e_1^1 \otimes \cdots \otimes e_{\alpha_{i_1}}^{i_1} \otimes \cdots \otimes e_{\alpha_{i_l}}^{i_l} \otimes \cdots \otimes e_1^d$, $1 \leq i_1<\cdots<i_l \leq d,  2 \leq \alpha_{i_t} \leq n_{i_{t}}, 1\leq t \leq l, 2\leq l \leq d$.

\medskip\noindent

Then $T_{p}\mathcal{C}={\rm span}\{p,E_{\alpha_{i}}^{i}: 1\leq i \leq d \}$.  By orthogonal decomposition, the normal space can be represented as:
\begin{equation}
   N_{p}\mathcal{C} =\operatorname{Span}\left\{E_{\alpha_{i_1}\cdots \alpha_{i_l}}^{i_{1}\cdots i_{l}}: 1 \leq i_1<\cdots<i_l \leq d,  2 \leq \alpha_{i_t} \leq n_{i_{t}}, 1\leq t \leq l, 2\leq l \leq d \right\}.
\end{equation}

\medskip\noindent

\textbf{3.2. The maximal value of $||A_{v}||^2$:} the second fundamental forms $B(E_{\alpha_{i}}^{i},E_{\beta_{i}}^{i})$ and $B(E_{\alpha_{i}}^{i},E_{\beta_{j}}^{j})$ can be computed as follows:

\medskip\noindent

We choose a curve $\gamma(t)$ starting from $p$ and has $E_{\alpha_{i}}^{i}$ as its tangent vector at $p$: 
$$
\gamma(t)=e_{1}^1 \otimes \cdots \otimes e_{1}^{i-1} \otimes \left({\rm cos} t \ e_{1}^{i}+ {\rm sin} t \ e_{\alpha_{i}}^{i}\right) \otimes e_{1}^{i+1}  \otimes \cdots \otimes e_1^{d}.
$$

Firstly, we consider an extension of $E_{\beta_{i}}^{i}$ along $\gamma(t)$: $\tilde{E}_{\beta_{i}}^{i}(\gamma(t))$, then derivative with $t$, we easily see that $\left( \frac{d}{dt}|_{t=0} \tilde{E}_{\beta_{i}}^{i}(\gamma(t))  \right)^{\perp}=0$, hence $B(E_{\alpha_{i}}^{i},E_{\beta_{i}}^{i})=0$.

\medskip\noindent

Secondly, we consider an extension of $E_{\beta_{j}}^{j}(j \neq i)$  along $\gamma(t)$: $\tilde{E}_{\beta_{j}}^{j}(\gamma(t))= e_1^1 \otimes \cdots \otimes \left({\rm cos} t \ e_{1}^{i}+ {\rm sin} t \ e_{\alpha_{i}}^{i}\right) \otimes \cdots \otimes e_{\beta_{j}}^{j} \otimes \cdots \otimes e_1^d $, then derivative in $t$, $ \frac{d}{dt}|_{t=0} \tilde{E}_{\beta_{j}}^{j}(\gamma(t)) =E_{\alpha_{i}\beta_{j}}^{ij}$ which is already located in $N_p \mathcal{C}$, hence $B(E_{\alpha_{i}}^{i},E_{\beta_{j}}^{j})=E_{\alpha_{i}\beta_{j}}^{ij}$, similarly, $B(E_{\beta_{j}}^{j},E_{\alpha_{i}}^{i})=E_{\alpha_{i}\beta_{j}}^{ij}$.

\medskip\noindent

We consider a unit normal $v \in N_p \mathcal{C}$: 
\begin{equation}\label{nv}
    v=\sum_{l\geq 2} a_{\alpha_{i_1}\cdots \alpha_{i_l}}^{i_1\cdots i_{l}}E_{\alpha_{i_1}\cdots \alpha_{i_l}}^{i_{1}\cdots i_{l}},
\end{equation}
then
$$
\langle B(E_{\beta_{j}}^{j},E_{\alpha_{i}}^{i}),v  \rangle=\langle B(E_{\alpha_{i}}^{i},E_{\beta_{j}}^{j}),v  \rangle=a_{\alpha_{i}\beta_{j}}^{ij},
$$
thus
$||A_{v}||^2=2\sum \left(a_{\alpha_{i}\beta_{j}}^{ij}\right)^2 \leq 2$, the equality is attained if and only if $a_{\alpha_{i_1}\cdots \alpha_{i_l}}^{i_1\cdots i_{l}}=0$ for $l\geq3$.

\medskip\noindent

\textbf{3.3. The normal radius:} as in Definition 2.1 and Lemma 2.2, we check at what time the tensor $T={\rm cos }t \ p+ {\rm sin} t \ v$ has rank $\textbf{R}(T)=1$, where $v$ is given as \eqref{nv}.

\medskip\noindent

\textbf{Claim:} The normal radius of Segre varieties are equal to $\frac{\pi}{2}$.

\medskip\noindent

\textbf{Proof of Claim:} The proof is easy, assume $\textbf{R}(T)=1$, say 
$$
T=(a_{1}e_{1}^{1}+\sum_{\alpha_{1}=1}^{n_{1}}b_{\alpha_{1}}^{1}e_{\alpha_{1}}^{1}) \otimes \cdots \otimes (a_{d}e_{1}^{d}+\sum_{\alpha_{d}=1}^{n_{d}}b_{\alpha_{d}}^{d}e_{\alpha_{d}}^{d}),
$$
we first assume there exists some $a_{i}=0$, then it implies ${\rm cos }t=0$, $t=\frac{\pi}{2}$, $v$ can be chosen to be those  $\textbf{R}(v)=1$.

\medskip\noindent

Then, assume all $a_{i}$ are not equal to zero, since $T$ contains no terms of the type $E_{\alpha_{i}}^{i}$, it implyies that
$$
\frac{a_{1}\cdots a_{d}}{a_i} \cdot b_{\alpha_{i}}^{i}=0 \ \ {\rm for \ \ any} \ \ 1\leq i \leq d,
$$
it concludes that all $b_{\alpha_{i}}^{i}=0$, $T=cp$ for some constant $c$, then must $c=-1$, $t=\pi$. $\Box$

\medskip\noindent

\textbf{3.4. Proof of the Main Theorem:} the regular part of Segre variety has a dimension $k={\rm dim } \ \mathcal{C}=\sum_{i=1}^{d} (n_i-1)+1$ from the parameterization given in \eqref{pa}, and $\alpha^2:={\rm sup}_{v}||A_{v}||^2=2$. We now estimate the vanishing angles by using Lawlor's table (cf. \cite[Sect. 1.4]{La91}) and formula \cite[Proposition 1.4.2]{La91}: 

(1): If $5 \leq  k \leq 12$, this rules out $C(\mathbb{S}^{1} \otimes \mathbb{S}^{1} \otimes \mathbb{S}^{1})$, then the vanishing angles are no more than $26.97^{\circ}$;

(2): If $k \geq 13$, let $\theta_{2}(k,\alpha)$ denote the estimated vanishing angle, $\theta_{2}(k,\alpha)$ is no less than the real vanishing angle and it decreases as $\alpha$ decreases (for example, see \cite[p. 376]{TZ20}), then \cite[Proposition 1.4.2]{La91} tells us:
$$
\tan \left(\theta_2(k, \sqrt{2})\right)<\frac{12}{k} \tan \left(\theta_2\left(12, \frac{12}{k}\sqrt{2}\right)\right)\leq \tan \left(\theta_2(12, \sqrt{2})\right)
$$
it implies that $\theta_2(k, \sqrt{2})\leq 8.08^{\circ}$ by Lawlor's table.

\medskip\noindent

Since the normal radius of those cones are $90^{\circ}$,  it follows from Theorem 2.6 that all Segre varieties are area-minimizing except for $C(\mathbb{S}^{1} \otimes \mathbb{S}^{1} \otimes \mathbb{S}^{1})$. $\Box$

\medskip\noindent
\medskip\noindent

\textbf{Acknowledgements:} This work is supported by the NSFC (Nos. 12301068, 12571057), the project of Stable Support for Youth Team in Basic Research Field, CAS (YSBR-001), and Xiaomi Young Scholar Fund. Part of the computations in the proof of Theorem 2.4 originates from an unpublished part of the author's PhD thesis at the University of Chinese Academy of Sciences. The author would like to express his sincere gratitude to his PhD supervisor, Prof. Xiaoxiang Jiao, for his invaluable guidance and constant encouragement.

	\medskip\noindent
		\medskip\noindent

\begin{flushleft}
			\medskip\noindent
			\begin{tabbing}
				XXXXXXXXXXXXXXXXXXXXXXXXXX*\=\kill
				Hongbin Cui\\
				School of Mathematical Sciences, University of Science and Technology of China\\
				Wu Wen-Tsun Key Laboratory of Mathematics, USTC, Chinese Academy of Sciences\\
				96 Jinzhai Road, Hefei, 230026, Anhui Province, China\\
				
				E-mail: cuihongbin@ustc.edu.cn
				
			\end{tabbing}
		\end{flushleft}        
        
  \end{document}